\documentclass[a4paper]{article}
\usepackage{amsfonts}
\usepackage{amsmath}
\usepackage{amssymb}
\usepackage{latexsym}
\usepackage{amstext}

\newcommand{\C}{\mathbb C}
\newcommand{\R}{\mathbb R}
\newcommand{\N}{\mathbb N}

\newtheorem{theorem}{Theorem}
\newtheorem{lemma}[theorem]{Lemma}
\newtheorem{prop}[theorem]{Proposition}
\newtheorem{cor}[theorem]{Corollary}
\newtheorem{definition}[theorem]{Definition}
\newtheorem{rem}[theorem]{Remark}

\newcommand{\mproof}{\noindent {\it Proof.~}}
\newcommand{\meproof}{\hfill {\blacksquare}}

\begin{document}
\title{On the Existence of the Fundamental Eigenvalue of an Elliptic Problem in $\R^N$ }
\author{J. Bellazzini, V. Benci M. Ghimenti, A.M.
Micheletti\\
{\it\small Dipartimento di Matematica Applicata} \\
{\it\small Universit\`a di Pisa}\\
{\it\small Via Buonarroti 1/C 56127 PISA -Italy}\\
{\it\small e-mails: j.bellazzini@ing.unipi.it,  benci@dma.unipi.it}\\
{\it\small e-mails: ghimenti@mail.dm.unipi.it,
a.micheletti@dma.unipi.it}}
\date{}
\maketitle

\begin{center}
{\bf\small Abstract}

\vspace{3mm}\hspace{.05in}\parbox{4.5in} {\small We study an
eigenvalue problem for functions in $\R^N$ and we find sufficient
conditions for the existence of the fundamental eigenvalue. This
result can be applied to the study of the orbital stability of the
standing waves of the nonlinear Schr\"odinger equation.}
\end{center}

\section{Introduction}

In this paper, we study the following eigenvalue problem:
\begin{equation}\tag{*}\label{eq}
\left\{
\begin{array}{ll}
 -\Delta u+F^{\prime }(u) =\lambda u  \\
\  u > 0
\end{array}
\right.
\end{equation}
\noindent where $u \in H^{1}(R^{n})$ and  $F:\R \rightarrow \R$ be
a even $C^{2}$ function such that $F(0)=F^{\prime }(0)=F^{\prime
\prime }(0)=0.$

In particular, we are interested in the existence of the
\emph{fundamental eigenvalue}; namely, the Lagrange multiplier of
the following minimization problem:

\

\noindent\textbf{Minimization problem}{\it: Find the minimum point
of the functional
\begin{equation*}
J(u)=\int_{^{N}}\frac{1}{2}|\nabla u|^{2}+F(u)dx
\end{equation*}
constrained to the manifold }
\begin{equation*}
M_{\rho }=\{v\in H^{1}(R^{N}),\ ||v||_{L^{2}}=\rho \}.
\end{equation*}

\bigskip

One of the motivations in studying the above problem is the
application to the following nonlinear Schr\"odinger equation:
\begin{equation}
i\frac{\partial \psi }{\partial t}+\Delta \psi -W^{\prime }(\left|
\psi \right| )\frac{\psi }{\left| \psi \right| }=0  \label{pina}
\end{equation}

\noindent where
\begin{equation*}
W(u)=\frac 12\Omega u^{2}+F(u).
\end{equation*}

We get the following result, which is a generalization of Theorem
II.2 of Cazenave and Lions \cite{CaLi}:

\begin{theorem}
\label{adone}{\it The nonlinear Schr\"odinger equation (\ref{pina})
admits orbitally stable standing waves if the above
minimization problem admits a solution with $\lambda<0$. These
standing waves have the form
\begin{equation*}
\psi \left( t,x\right) =u(x)e^{-i\omega t}
\end{equation*}
where $\omega =\lambda - \Omega $, and where $u$ and $\lambda $ are
the first eigenfunction and the first eigenvalue of problem
(\ref{eq}) respectively.}
\end{theorem}

We notice that if $\lambda>0$,  it has been proved in \cite{BeLi}
that (\ref{eq}) has no radial solution. Clearly the above result has
a limited interest if we do not know when the fundamental solution
of (\ref{eq}) exists.

Now we will state the main result of this paper concerning problem
(\ref{eq}). We state the hypothesis
\begin{equation}
|F^{\prime }(s)|\leq c_{1}|s|^{q-1}+c_{2}|s|^{p-1}\text{ for some
}2<q\leq p<2^{\ast }.  \tag{$F_{p}$}  \label{Fp}
\end{equation}
We also assume
\begin{equation}
F(s)\geq -c_{1}s^{2}-c_{2}|s|^{\gamma }\text{ for some
}c_{1},c_{2}\geq 0,\ { \gamma <2+\frac{4}{N}}  \tag{$F_{0}$}
\label{F0}
\end{equation}
and
\begin{equation}
\text{there exists }s_{0}\in \text{ such that }F(s_{0})<0.
\tag{$F_{1}$} \label{F1}
\end{equation}

\begin{theorem}
\label{mainradiale} {\it Let $F$ satisfy (\ref{Fp}), (\ref{F0}) and
(\ref{F1}). Then, $\exists\ \bar \rho$ such that $\forall\
\rho>\bar\rho$ there exists $ \bar u\in H^1$ satisfying

\begin{equation*}
J(\bar u)= \inf_{\{v\in H^1,\ ||v||_{L^2}=\rho\}}J(v),
\end{equation*}
with $||\bar u||_{L^2}=\rho$. Then, there exist $\lambda$ and $\bar
u $ that solve (\ref{eq}), with $ \lambda<0$ and $\bar u $ positive
radially symmetric.}
\end{theorem}

In order to have stronger results, we can replace (\ref{F1}) with
the following hypothesis

\begin{equation}  \label{F2}
F(s)<-s^{2+\epsilon},\ 0<\epsilon<\frac 4N\text{ for small }s.
\tag{$F_2$}
\end{equation}

In this case we find the following results concerning the
existence of the minimizer of $J(u)$ for any $\rho$.

\begin{cor}
\label{cor1} {\it If (\ref{Fp}), (\ref{F0}) and (\ref{F2}) hold,
then for all $ \rho$, there exists $\bar u\in H^1$, with $||\bar
u||_{L^2}=\rho$, such that
\begin{equation*}
J(\bar u)= \inf_{\{v\in H^1, ||v||_{L^2}=\rho\}}J(v).
\end{equation*}}
\end{cor}

In particular, for $N=3$ we have

\begin{cor}
\label{cor2} {\it Let $N=3$. If (\ref{Fp}) and (\ref{F0}) hold and
$F\in C^{3}$, with $F^{\prime \prime \prime }(0)<0$, then for all
$\rho $, there exists $\bar{u} \in H^{1}$ with
$||\bar{u}||_{L^{2}}=\rho $ such that
\begin{equation*}
J(\bar{u})=\inf_{\{v\in H^{1},||v||_{L^{2}}=\rho \}}J(v).
\end{equation*}}
\end{cor}

There is a lot of literature for nonlinear eigenvalue problems. We
refer to Rabinowitz \cite{Rab} and the references therein. However,
as far as we know, there are no results in the case when the
nonlinearity is not a compact perturbation of the Laplace operator.

In \cite{BeLi}, in order to prove the existence of a solution for
the problem (\ref{eq}) with $\lambda$ negative and fixed, the
authors used a slightly weaker version of (\ref{Fp}) and a slightly
different version of (\ref{F1}). In fact, (\ref{F0}) is used in
\cite{CaLi} in order to obtain that the Cauchy problem associated to
equation (\ref{pina}) has a solution for all time, and also to prove
the orbital stability of the standing wave relative to the ground
state solution.

The paper is organized as follows. In section \ref{P} we prove
theorem \ref {mainradiale}. In section \ref{A} we prove theorem
\ref{adone}. This theorem is a generalization of theorem II.2 of
\cite{CaLi}. As a matter of fact, theorem \ref{adone} can be
obtained following the same type of argument as in Theorem II.2, as
claimed in \cite{CaLi}. Here we give a complete and different proof
which is based on the ''splitting lemma''. In the appendix we prove
the splitting lemma in the form used here.

\section{Proofs of the main results\label{P}}

\begin{lemma}\label{lemmone}
{\it If $F$ satisfies (\ref{F1}), then $\exists\  \bar \rho$ such
that $\forall\ \rho>\bar\rho$
$$\inf\limits_{||u||_{L^2}=\rho}J(u)<0.$$
Otherwise, if $F$ satisfies (\ref{F2}) then
$\inf\limits_{||u||_{L^2}=\rho}J(u)<0$ for all $\rho$.}
\end{lemma}

\mproof We build a sequence of radial functions $u_n$ in $H^1$
such that $J(u_n)<0$ for large $n$. The sequence is defined as
follows:
\begin{equation}
u_n(r)= \left\{
\begin{array}{lll}
s_0&&r<R_n;\\
s_0-s_0(r-R_n)&&R_n\leq r\leq R_n+1;\\
0&& r>R_n+1.
\end{array}
\right.
\end{equation}

We show that $J(u_n)<0$ when $R_n \rightarrow +\infty$. Notice that
$|\nabla u_n|^2=|\frac{\partial u_n}{\partial r}|^2=s_0^2$. We have
$$\int_{\R^N} \frac12 |\nabla u_n|^2+F(u_n)dx $$ $$\leq   C_1
\int_{R_n}^{R_n+1} \left[ \left|\frac{\partial u_n}{\partial
r}\right|^2+ \sup_{|s|\in[R_n,R_n+1]} F(s)\right] r^{N-1}dr  +
C_2\int_0^{R_n}F(s_0)r^{N-1}dr \vspace{2mm}\noindent$$
 with $C_1$ and $C_2$ strictly positive. This
proves the first statement.

Now, we want to prove that if (\ref{F2}) holds, then
$$\inf\limits_{||u||_{L^2}=\rho}J(u)<0$$ for all $\rho$. We use the
same approach as before; we build a sequence of radial functions
that are constant in a ball with a suitable cut-off. Let $u_n$ be
\begin{equation}
u_n(r)= \left\{
\begin{array}{lll}
s_n&&r<R_n;\\
s_n-\frac{s_n}{R_n}(r-R_n) &&R_n\leq r \leq 2R_n;\\
0&& r>2R_n.
\end{array}
\right.
\end{equation}

We study $J(u_n)$ when $R_n \rightarrow +\infty$; and, due to the
constraint $||u_n||_{L^2}=\rho$, we have
\begin{equation}\label{mr}
\lim_{R_n \rightarrow \infty} s_n^2 R_n^N = \gamma > 0.
\end{equation}
We can choose $R_n$ sufficiently large such that $F(u_n)\leq 0$.
Therefore,
\begin{eqnarray*}
J(u_n)&=&\int_{\R^N} \frac12 |\nabla u_n|^2+F(u_n)dx = \\
&=& C \int_{R_n}^{2R_n} \left|\frac{\partial u}{\partial
r}\right|^2r^{N-1} dr +
C \int_{0}^{R_n} F(u_n)r^{N-1}dr +\\&+&C \int_{R_n}^{2R_n} F(u_n)r^{N-1}dr \leq\\
&\leq& C \int_{R_n}^{2R_n} \left|\frac{\partial u}{\partial
r}\right|^2r^{N-1} dr +
C \int_{0}^{R_n} F(u_n)r^{N-1}dr =\\
&=& C\int_{R_n}^{2R_n} \frac{s_n^2}{R_n^2}r^{N-1}dr +
C\int_{0}^{R_n} F(s_n)r^{N-1}dr \leq\\
&\leq& C  s_n^2R_n^{N-2}
-C \int_{0}^{R_n} s^{2+\epsilon}r^{N-1}dr=\\
&=& Cs_n^{2}R_n^{N-2}-Cs^{2+\epsilon}R_n^{N},
\end{eqnarray*}
where $C$ are positive constants. By (\ref{mr}) we have
$R_n=O(s_n^{-2/N})$ and thus
\begin{equation}
J(u_n)\leq O(s_n^{4/N})-O(s_n^\epsilon)\rightarrow0^-,
\end{equation}
if $0<\epsilon<4/N$. $\meproof$

\begin{rem}
In the proof of previous lemma we have used radially symmetric
functions. So, in the same way we can obtain that
\begin{equation}
\inf\limits_{\{u\in H^1_r, ||u||_{L^2}=\rho\}}J(u)<0
\end{equation}
with the same hypothesis.
\end{rem}

\begin{prop}\label{max}
{\it If there exist $c_1,c_2 \geq 0$ such that
\begin{equation}\tag{$F_0$}
F(u)\geq -c_1u^2-c_2u^{\gamma}, \ 2\leq\gamma<2+\frac{4}{N}
\end{equation}
then
\begin{enumerate}
\item $\inf\limits_{||u||_{L^2}=\rho}J(u)>-\infty$, and \item any
minimizing sequence $u_n$, i.e. $J(u_n)\rightarrow c$,
$||u_n||_{L^2}=\rho$, is bounded in $H^1$.
\end{enumerate}}
\end{prop}

\mproof We apply the Sobolev inequality (see\cite{VoPa})
\begin{equation}\label{Sobolev}
||u||_{L^q} \leq
b_q||u||_{L^2}^{1-\frac{N}{2}+\frac{N}{q}}||\nabla
u||_{L^2}^{\frac{N}{2}-\frac{N}{q}}
\end{equation}
that holds for $2 \leq q \leq 2^{*}$ when $N \geq 3$.

 From equation
(\ref{Sobolev}) we have that any function $u$ such that
$||u||_{L^2}=\rho$ fulfills the following equation:
\begin{equation}\label{sopra}
||u||_{L^q}^q \leq b_{q,\rho} ||\nabla u||_{L^2}^{\frac{qN}{2}-N}.
\end{equation}
Now we notice that, by (\ref{sopra}), for all $u\in H^1$ with
$||u||_{L^2}=\rho$
\begin{eqnarray}
  J(u)&\geq& \int \frac 12 |\nabla u|^2 -c_1u^2-c_2u^{\gamma}dx\\
&\geq &\int \frac 12 |\nabla u|^2 dx- c_2
b_{\gamma,\rho}\left(\int{|\nabla u|^2}dx\right)^{\frac{\gamma
N}{2}-N} -c_1\rho^2\nonumber.
\end{eqnarray}
If $\frac{\gamma N}{2}-N<2$, i.e if $\gamma < 2+\frac{4}{N}$, we
have
\begin{equation}
 J(u)\geq \frac{1}{2}||\nabla u||_{L^2}^2 +O(||\nabla u ||_{L^2}^2).
\end{equation}
The proof follows easily. $\meproof$

\begin{lemma}\label{basiclemma}
{\it Let (\ref{Fp}) hold, and let $\bar{u} \neq 0$, $\bar{u}\in
H^1(\R^N)$ and $\bar{\lambda}\in \R$ be such that
\begin{equation}\tag{*}
-\Delta \bar{u}+F'(\bar{u})=\bar{\lambda}\bar{u}.
\end{equation}
If
\begin{equation}\label{Ju}
J(\bar{u})=\int_{\R^N} \frac 12|\nabla \bar{u}|^2+F(\bar{u})dx <
0,
\end{equation}
then
$$\bar{\lambda}<0.$$}
\end{lemma}

\mproof By assumption (\ref{Fp}), if $\bar u \in H^1$ solves
(\ref{eq}), by a bootstrap argument, $\bar u\in
H^{2,2}_{loc}(\R^N)$. Furthermore, since $F(u)\in L^1$, we can apply
the Derrick-Pohozaev identity (see \cite{De,Po,BeLi})
\begin{equation}\label{pohozaev}
\int_{\R^N} |\nabla \bar{u}|^2=
\frac{2N}{N-2}\int_{\R^N}\frac{\bar{\lambda}\bar{u}^2}{2}-F(\bar{u})
dx.
\end{equation}
The function $\bar{u}$ satisfies the equation $-\Delta
\bar{u}+F'(\bar{u})=\bar{\lambda}\bar{u}$. Therefore, by
integration, we get
\begin{equation}\label{bara}
\frac 12 \int_{\R^N} |\nabla \bar{u}|^2 dx = \frac 12 \int_{\R^N}
\bar{\lambda}\bar{u}^2-F'(\bar{u})\bar{u}dx.
\end{equation}
By equation (\ref{pohozaev}) we have
\begin{equation}\label{pohozaev2}
\frac{N-2}{2N} \int_{\R^N} |\nabla \bar{u}|^2 dx =
\int_{\R^N}\frac{\bar{\lambda}\bar{u}^2}{2}-F(\bar{u}) dx,
\end{equation}
and subtracting (\ref{bara}) from (\ref{pohozaev2}), we get
\begin{equation}
-\frac{1}{N}\int_{\R^N} |\nabla \bar{u}|^2 =
\left[\frac{N-2}{2N}-\frac 12 \right]\int_{\R^N} |\nabla
\bar{u}|^2 = \int_{\R^N} \frac 12 F'(\bar{u})\bar{u}- F(\bar{u}).
\end{equation}
This proves that
\begin{equation}\label{derpoho}
\int_{\R^N} \frac 12 F'(\bar{u})\bar{u}- F(\bar{u})dx <0.
\end{equation}
On the other hand, by (\ref{Ju}) and (\ref{bara})
\begin{equation}
2J(\bar{u})+ \int_{\R^N} F'(\bar{u})\bar{u}- 2F(\bar{u})dx =
\bar{\lambda}  \int_{\R^N}\bar{u}^2dx,
\end{equation}
and thus we have $\bar{\lambda} <0$. $\meproof$

\begin{rem}
With the same argument as before, we can prove that
\begin{equation}\label{amra}
\frac{2}{N-2}\bar{\lambda}\int_{\R^N}\bar{u}^2dx= \int_{\R^N} 2^*
F(\bar{u}) -F'(\bar{u})\bar{u}dx
\end{equation}
without any assumption on the sign of $J(\bar{u})$. Indeed, we have,
as above,
\begin{eqnarray}
\int_{\R^N} |\nabla \bar{u}|^2 dx &=& \frac{N}{N-2}\bar{\lambda}
\int_{\R^N}\bar{u}^2dx -
\frac{2N}{N-2}\int_{\R^N}F(\bar{u})dx\\
&&\nonumber\\
 \int_{\R^N} |\nabla \bar{u}|^2 dx &= &
   \int_{\R^N} \bar{\lambda}\bar{u}^2-F'(\bar{u})\bar{u}dx.
\end{eqnarray}
By subtraction we obtain (\ref{amra}). Notice that this is an
Ambrosetti-Rabinowitz type inequality.
\end{rem}

\begin{rem}\label{remmino}
We notice that the infimum of $J(u)$ for  $||u||_{L^2}=\rho, u\in
H^1(\R^N)$ or $||u||_{L^2}=\rho, u\in H^1_r(\R^N)$ has the same
value. If we call
\begin{eqnarray*}
&&c=\inf \{J(u), \ ||u||_{L^2}=\rho, u\in H^1(\R^N)\}\\
&&c_r=\inf \{J(u), \ ||u||_{L^2}=\rho, u\in H^1_r(\R^N)\}
\end{eqnarray*}
 we have $c=c_r$. Clearly
$c\leq c_r$ and, if $u_n \geq 0$ is a minimizing sequence for $c$,
we denote by $u_n^*$ the Schwartz spherical rearrangement of
$u_n$. Now, $u_n^* \in H^1_r(\R^N)$, $||u_n^*||_{L^2}=\rho$, and
$J(u_n^*)\leq J(u_n)$, thus $u_n^*$ is a minimizing sequence.
Therefore, $c=c_r$.
\end{rem}

\begin{prop}\label{lambda}
{\it Let (\ref{F0}) and (\ref{Fp}) hold; and Let $u_n$  be a
minimizing P-S sequence such that $J(u_n)\rightarrow c$, where
$c=\inf_{||u||_{L^2}=\rho}J(u)<0$. Then there exists a sequence
$\lambda_n$ of Lagrange multipliers such that
\begin{equation}\label{ps}
-\Delta u_n +F'(u_n)-\lambda_nu_n=\sigma_n \rightarrow  0.
\end{equation}
We have that $\lambda_n$ is bounded in $\R$}.
\end{prop}

\mproof By (\ref{ps}) and Proposition \ref{max}, we have
\begin{equation*}
\left| \int_{\R^N} |\nabla u_n|^2 +F'(u_n)u_n-\lambda_n
u_n^2dx\right| \leq ||\sigma_n||_{H^{-1}}||u_n||_{H^1}\rightarrow
0.
\end{equation*}
We have
\begin{eqnarray*}
&&\int_{\R^N} |\nabla u_n|^2 +F'(u_n)u_n-\lambda_n u_n^2dx=\\
&&\int_{\R^N} |\nabla u_n|^2 +2F(u_n)-2F(u_n)+F'(u_n)u_n-\lambda_n u_n^2dx=\\
&&2J(u_n)-\lambda_n\rho^2+\int_{\R^N} F'(u_n)u_n-2F(u_n)dx
\rightarrow0.
\end{eqnarray*}
Furthermore, we have that $J(u_n)$ is bounded; and also by
(\ref{Fp}),
$$\left|\int_{\R^N} F'(u_n)u_n-2F(u_n)dx\right|
\leq C\left(||u_n||_{H^1}^q+||u_n||_{H^1}^p\right)<+\infty. $$ Then
$\lambda_n$ is bounded and statement is proved. $\meproof$

\vspace{3mm} \noindent{\it Proof of Theorem \ref{mainradiale}}.  For
the Palais principle the critical point of the functional on the
manifold $\{||u||_{L^2}=\rho$, $u\in H^1_r\}$ are still critical
points on the manifold $\{||u||_{L^2}=\rho$, $u\in H^1\}$. Remark
\ref{remmino} assures that the minimizers of the functional on the
manifold $\{||u||_{L^2}=\rho$, $u\in H^1_r\}$ are still minimizers
of the functional in $H^1$. So, we study the existence of minimizers
in $H^1_r$.

Let $u_n$  be a minimizing sequence such that $||u_n||_{L^2}=\rho$;
$F(s)$ is an even function, we can even take $u_n \geq 0$. By Lemma
\ref{lemmone}, we can take $\rho$ sufficiently large such that
$J(u_n) \rightarrow c<0$. By the Ekeland principle, we can assume
that $u_n$ is a Palais-Smale sequence for the functional $J$
restricted on the manifold $||u_n||_{L^2}=\rho$.

By Lemma \ref{max} and Proposition \ref{lambda}, $u_n$ is bounded in
$H^1_r$ and $\lambda_n$ is bounded in $\R$. We have
\begin{eqnarray*}
&& \lambda_n \rightarrow \bar{\lambda}\\
&&  u_n \rightharpoonup u \ \ \ \text{weakly in}\ H^1_r(\R^N)\\
&& u_n \rightarrow u \ \ \  \text{strongly in }\  L^p(\R^N), \ 2<p<2^* \\
&&  u_n \rightarrow u \ \ \ \text{strongly in }\  L^p(B), \ B\text{
compact} \ 2 \leq p < 2^*.
\end{eqnarray*}
Moreover, for any radially symmetric function we have the following
decay when $|x|\rightarrow \infty$:
\begin{equation}\label{decadenz}
|u_n(x)|\leq \alpha\frac{||u_n||_{H^1}}{|x|^{\frac{N-1}{2}}}\ \
\text{for} \ |x|>\beta
\end{equation}
where $\alpha, \beta$ depend only on $N$ (see  for instance \cite{BeLi}).\\

By (\ref{Fp}), it is easy to see that
\begin{equation}
-\Delta u +F'(u)=\bar{\lambda}u.
\end{equation}
Indeed, we have for any $\varphi \in C^{\infty}_0(\R^N)$ radially
symmetric
\begin{equation}
\int_{\R^N} \nabla u_n \nabla \varphi dx
-\int_{\R^N}F'(u_n)\varphi dx -\lambda_n\int_{\R^N} u_n\varphi dx
\rightarrow 0
\end{equation}
as $n \rightarrow +\infty$. By (\ref{Fp}), we have that
$$\int_{\R^N} F'(u_n)\varphi\rightarrow \int_{\R^N} F'(u)\varphi. $$
Then, as $n \rightarrow +\infty$,
\begin{eqnarray}
&&\int_{\R^N} \nabla u_n \nabla \varphi dx
-\int_{\R^N}F'(u_n)\varphi dx -\lambda_n\int_{\R^N} u_n\varphi dx
\rightarrow  \\ \nonumber && \rightarrow \int_{\R^N} \nabla u
\nabla \varphi dx -\int_{\R^N}F'(u)\varphi dx
-\bar{\lambda}\int_{\R^N} u\varphi dx.
\end{eqnarray}
We show that $u \neq 0$.

Indeed, the Neminski operator
$$F: L^t(\R^N) \rightarrow L^1(\R^N)  \ 2<t<2^*$$
is continuous by (\ref{Fp}) and $u_n \rightarrow u$ in $L^t(\R^N),
\ 2<t<2^*$. Hence we have
\begin{equation}
\frac{1}{2} \int_{\R^N} |\nabla u|^2+\int_{\R^N} F(u)dx \leq
\lim_{n\rightarrow \infty}J(u_n)=c<0,
\end{equation}
which proves that $u \neq 0$.

At this point, we have that $u\neq 0$ is a weak solution of
\begin{equation}
-\Delta u +F'(u)=\bar \lambda u
\end{equation}
and that
\begin{equation}
J(u)<0.
\end{equation}
Thus the hypotheses of Lemma \ref{basiclemma} are fulfilled, and
$\bar \lambda<0$.

Considering two functions $u_n$ and $u_m$ in the minimizing P-S
sequence, we have
\begin{eqnarray*}
&& -\Delta u_n +F'(u_n)-\lambda_nu_n=\sigma_n \rightarrow  0\\
&& -\Delta u_m +F'(u_m)-\lambda_mu_m=\sigma_m \rightarrow  0.
\end{eqnarray*}
By subtraction we get
\begin{equation}
 -\Delta (u_n-u_m) +F'(u_n)-F'(u_m)-\bar{\lambda}(u_n-u_m) \rightarrow  0,
\end{equation}
and we obtain
\begin{equation}\label{converg}
\int_{\R^N}
|\nabla(u_n-u_m)|^2dx+\int_{\R^N}(F'(u_n)-F'(u_m))(u_n-u_m)-
\bar{\lambda}(u_n-u_m)^2dx \rightarrow 0.
\end{equation}
On any compact ball $B$ we have, by standard arguments, that
\begin{eqnarray*}
&&\int_{B}(F'(u_n)-F'(u_m))(u_n-u_m)\rightarrow0\\
&&\int_{B}\bar{\lambda}(u_n-u_m)^2dx \rightarrow 0.
\end{eqnarray*}
Then,
\begin{eqnarray}\label{converg2}
\int_{\R^N} |\nabla(u_n-u_m)|^2dx&+&\int_{B^c}(F'(u_n)-F'(u_m))(u_n-u_m)-\\
&-&\int_{B^c}\bar{\lambda}(u_n-u_m)^2dx\rightarrow0.\nonumber
\end{eqnarray}
By lemma \ref{basiclemma}, we have  $\bar{\lambda}<0$ and
\begin{eqnarray}\nonumber
&&\int_{B^c}(F'(u_n)-F'(u_m))(u_n-u_m)-\bar{\lambda}(u_n-u_m)^2dx=\\
&=&\int_{B^c} (F''(\theta u_n+(1-\theta)
u_m)-\bar{\lambda})(u_n-u_m)^2dx.
\end{eqnarray}
Thus, remembering that $F''(0)=0$ and due to (\ref{decadenz}), we
have $F''(\theta u_n+(1-\theta)u_m)<<1$, and
\begin{equation}
F''(\theta u_n+(1-\theta) u_m)-\bar{\lambda}>0
\end{equation}
for $B$ sufficiently large. By equation (\ref{converg}) we get
\begin{equation}
\int_{\R^N} |\nabla(u_n-u_m)|^2dx\rightarrow 0
\end{equation}
\begin{equation}
\int_{\R^N} |(u_n-u_m)|^2dx\rightarrow 0.
\end{equation}
Then the sequence $\{u_n\}_n$ is a Cauchy sequence in $H^1_r(\R^N)$;
thus $u_n\rightarrow u$ strongly in $H^1(\R^N)$ and
$||u||_{L^2}=\rho$.

$\meproof$

The proofs of Corollary \ref{cor1} and Corollary \ref{cor2} are
straightforward. Moreover, we prove a non existence result when
$F(s)=s^p$ with $2<p<2^*$.
\begin{rem}
Let $F$ satisfy (\ref{Fp}) and (\ref{F0}). If
\begin{equation}\label{nonex}
0 \leq 2F(s) \leq  F'(s)s \ \text{for all} \ s
\end{equation}
then (\ref{eq}) has no nontrivial solution in $H^1(\R^N)$ for all
$\lambda$. The proof is a consequence of the Derrick-Pohozaev
identity.

 Let us suppose that there exists $u \in H^1(\R^N)$, $u \neq 0$,
and $\lambda$ such that (\ref{eq}) holds: by bootstrap arguments we
have that $u\in H^{2,q}$ for all $q$ and we can apply the
Derrick-Pohozaev identity. Therefore, by (\ref{derpoho}) no solution
of (\ref{eq}) can  satisfy (\ref{nonex}).
\end{rem}

\section{Stability of the nonlinear Schr\"odinger equation \label{A}}

We consider the nonlinear Schr\"odinger equation

\begin{equation}\label{schr}\tag{$\dagger$}
\left\{
\begin{array}{ll}
i\frac{\partial \psi}{\partial t}+\Delta\psi-F'(\psi)=0&(t,x)\in \R^+\times\R^N\\
&\\
\psi(0,x)=\psi_0(x)&
\end{array}
\right.
\end{equation}
where $F:\C\rightarrow\R$ is a an even radial function such that
$F(|\xi|)$ satisfies (\ref{F0}),(\ref{Fp}) and (\ref{F1}). It is
well known that there exists a unique solution $\psi \in
C([0,+\infty),H^1(\R^N))$, see \cite{Str78,Gin79, Lin78}. We notice
that problem (\ref{pina}) reduces to (\ref{schr}) in the case
$\Omega=0$. It is easy to see that this hypothesis is not
restrictive. Setting $\psi=u(t,x)e^{iS(t,x)}$ we have that any
solution of (\ref{schr}) verifies
\begin{equation}\label{continuity}
\left\{
\begin{array}{l}
\frac 12 \frac{d}{dt}u^2 + \nabla \cdot (u^2 \nabla S)=0\\
\\
u\partial_t S-\Delta u + u|\nabla S|^2+F'(u)=0.
\end{array}
\right.
\end{equation}
It is well known that (\ref{continuity}) are the Euler-Lagrange
equations of the action functional given by
\begin{equation}
A(u,S)=\iint \frac{1}{2}u^2 \partial_t S + \frac 12 |\nabla u|^2+
\frac 12 u^2 |\nabla S|^2+F(u) dtdx.
\end{equation}
Since the energy is given by
\begin{equation}
E(\psi)=E(u,S)=\int\frac 12 |\nabla u|^2+ \frac 12 u^2 |\nabla
S|^2+F(u) dx,
\end{equation}
any solution of (\ref{continuity}) satisfies
\begin{eqnarray}
&&\frac{d}{dt} \int u(t,x)^2dx=0\label{normcost}\\
&&\frac{d}{dt}E(u,S)=0.\label{energcost}
\end{eqnarray}
Hence, for any solution of (\ref{schr}), equations (\ref{normcost})
and (\ref{energcost}) become
\begin{eqnarray}
&&||\psi(t,\cdot)||_{L^2(\R^N)}=||\psi_0(\cdot)||_{L^2(\R^N)}\\
&&E(\psi(t,x))=E(\psi_0)
\end{eqnarray}
for all $t$.

A  solution of (\ref{schr})  is called stationary solution if
$\psi=v(x)e^{-i\omega t}$. Such a solution satisfies the nonlinear
eigenvalue problem
\begin{equation}\tag{*}
-\Delta v +F'(v)=\omega v.
\end{equation}
By Theorem \ref{mainradiale}, if $F$ satisfies (\ref{F0}),
(\ref{F1}) and (\ref{Fp}), there exist $(\bar{u},\lambda)$ that
satisfies (\ref{eq}) such that
\begin{equation*}
J(\bar{u})= \inf_{\{v\in H^1,\ ||v||_{L^2}=\rho\}}J(v),
\end{equation*}
with $||\bar{u}||_{L^2}=\rho$, for some $\rho$, and $\lambda$ is a
Lagrange multiplier. So we have that $\psi=\bar{u}(x)e^{-i\lambda
t}$ is a stationary solution of (\ref{schr}) with initial condition
$\psi(0,x)=\bar{u}(x)$. Notice that for stationary solution,
$E(\psi)=J(u)$. Indeed, we have
\begin{equation}
E(u(t,x)e^{iS(t,x)})=J(u(t,x))+\frac 12 \int u^2|\nabla S|^2dx
\end{equation}
for all $t$. Now we prove the orbital stability of the stationary
solution found in the previous section.

We define

\begin{equation}
S=\left\{u(x)e^{i\theta};\ \theta \in S^1, ||u||_{L^2}=\rho, J(u)=
\inf_{\{v\in H^1,\ ||v||_{L^2}=\rho\}}J(v) \right\}.
\end{equation}
Clearly, for any $q\in \R^N$ we have that $\bar u(x+q)\in S$.

\begin{definition}
$S$ is {\it orbitally stable} if
\begin{equation*}
\forall \varepsilon, \ \exists \delta>0 \ \text{s.t.} \  \forall
\psi_0 \in H^1(\R^N), \ \inf_{u\in S} ||\psi_0-u||_{H^1}<\delta \
\text{implies}
\end{equation*}
\begin{equation*}
\forall t\geq 0 \ \ \inf_{u\in S}
||\psi(t,\cdot)-u||_{H^1}<\varepsilon
\end{equation*}
where $\psi(t,x)$ is the solution of (\ref{schr}) with initial
data $\psi_0$.
\end{definition}
Let us suppose that  $S$ is not orbitally stable, i.e that
\begin{equation*}
\exists \varepsilon, \ \exists  \psi_n(0,x) \in H^1(\R^N), \
\inf_{u\in S} ||\psi_n(0,x)-u||_{H^1}\rightarrow 0 \
\text{implies}
\end{equation*}
\begin{equation*}
\exists t_n \geq0 \ \ \inf_{u\in S}
||\psi_n(t,\cdot)-u||_{H^1}>\varepsilon.
\end{equation*}
We can suppose that
$||\psi_n(t_n,x)||_{L^2}=||\psi_n(0,x)||_{L^2}=\rho$. Indeed, if
$$||u_n(t,x)||_{L^2} \rightarrow \rho$$ there exists a sequence
$\alpha_n=\frac{\rho}{||u_n||_{L^2}}$ such that $||\alpha_n
u_n||_{L^2}=\rho$ and $J(\alpha_nu_n) -J(u_n) \rightarrow 0$. We
notice that, denoting $\psi_n(t,x)=u_n(t,x)e^{iS_n(t,x)}$,
\begin{equation}
J(u_n(0,x))\rightarrow J(\bar{u}),
\end{equation}
i.e $u_n(0,x)$ is a minimizing sequence of $J(u)$ on
$||u||_{L^2}=\rho$. Moreover,
\begin{equation}
E(\psi_n(t_n,x))=E(\psi_n(0,x))\rightarrow
E(\bar{u}(0,x))=J(\bar{u}).
\end{equation}
Hence we have that $$u_n(t_n,x)\ {\it is\,\, a \,\,minimizing\,\,
sequence\,\, on\,\, } ||u||_{L^2}=\rho.$$

Now we prove that any minimizing sequence for $J(u)$ on
$||u||_{L^2}=\rho$ does converge in $H^1$. This proves clearly that
$S$ is orbitally stable. As a matter of fact this result can be
proved, as claimed in \cite{CaLi}, as a consequence of the
concentration-compacteness principle of P.L. Lions
\cite{Li01}-\cite{Li02}. Here, in order to give a self contained and
simpler formulation, we prove a ``Splitting Lemma" which describes
the behaviour the Palais-Smale sequences. This lemma is a well known
result of Struwe \cite{Str84}. To prove this Lemma we make the
following growth assumption
\begin{equation}\label{Fp'}\tag{$F_p'$}
|F''(s)|\leq c_1|s|^{q-2}+c_2|s|^{p-2}\text{ for some }2<q\leq
p<2^*.
\end{equation}
We know that every critical point of $J$ on $||u||_{L^2}=\rho$ is
still a critical point of a corresponding functional
\begin{equation}
J_{\lambda}(u)=J(u)-\lambda\int_{\R^N}u^2 dx
\end{equation}
where $\lambda$ is the suitable Lagrange multiplier.
\begin{lemma}\label{splitting}
{\it Let (\ref{F0}) and (\ref{Fp'}) hold, and let $u_n$ be a
Palais-Smale sequence for $J_{\lambda}$ with $\lambda < 0$ and
$||u_n||_{L^2}\rightarrow \rho$. Then there exist $k$ sequence of
points $\{y^j_n\}_{n\in \N}$ $(1\leq j\leq k)$ with
$|y^j_n|\rightarrow +\infty$ such that, up to a subsequence:
\begin{enumerate}
\item $u_n=u^0+\sum_j u^j(x+y^j_n)+w_n$ with $w_n \rightarrow 0$
in $H^1$ \item $||u_n||_{L^2}^2 \rightarrow ||u^0||^2_{L^2}+\sum_j
||u^j||^2_{L^2}$ \item $J_{\lambda}(u_n) \rightarrow
J_{\lambda}(u^0)+\sum_j J_{\lambda}(u^j)$
\end{enumerate}
where $u^0$ and $u^j$ are weak solutions of (\ref{eq}).}
\end{lemma}

The Proof of Lemma \ref{splitting} is given in Appendix.
\begin{prop}\label{lions}
{\it Suppose that for any $\rho$ there exists
$I_{\rho^2}:=\min\limits_{||u||^2_{L^2}=\rho^2}J(u)$. Then, for any
$\mu \in (0,\rho)$ we have
\begin{equation}\label{idealions}
I_{\rho^2}<I_{\mu^2}+I_{\rho^2-\mu^2}.
\end{equation}}
\end{prop}
\mproof We prove that $I_{\theta\rho^2}<\theta I_{\rho^2}$ for any
$\rho>0$ and for any $\theta>1$. We take $u$ such that
$J(u)=I_{\rho^2}$. Thus,
$||u(\frac{x}{\theta^{1/N}})||^2_{L^2}=\theta \rho^2$. We have
\begin{eqnarray}
I_{\theta \rho^2}\leq
J\left(u\left(\frac{x}{\theta^{1/N}}\right)\right)=
\theta \left( \int_{R^N} \frac{1}{2}\left(\frac{1}{\theta}\right)^{2/N} |\nabla u|^2 + F(u) dx\right)\\
 <  \theta \left(\int_{\R^N} \frac12  |\nabla u|^2+ F(u)dx \right)=\theta I_{\rho^2}. \nonumber
\end{eqnarray}
By simple arguments we obtain (\ref{idealions}). In fact, if
$h(x)$ is a real function such that, for all $x>0$ and for all
$\theta>1$
\begin{equation}
h(\theta x)<\theta h(x)
\end{equation}
then we have, for all $y\in(0,x)$
\begin{equation}
h(x)<h(y)+h(x-y).
\end{equation}
$\meproof$

With  Lemma \ref{splitting} and Proposition \ref{lions} we can prove
the following Theorem.
\begin{theorem}\label{stability}
{\it Let (\ref{F0}), (\ref{F2}) and (\ref{Fp'}) hold. Then $S$ is
orbitally stable.}
\end{theorem}
\mproof Let, as before,  $u_n(t_n,x)$ be a minimizing sequence of
$J$ on $||u||_{L^2}=\rho$. By the Ekeland principle we can suppose
that it is a Palais-Smale sequence for $J$ and, thus, a Palais-Smale
sequence for $J_\lambda$ with $\lambda<0$, as proved in Theorem
\ref{mainradiale}.

By Lemma \ref{splitting} and Proposition \ref{lions} we have two
cases:
\begin{enumerate}
\item $u_n=u^0+w_n$  with $w_n \rightarrow 0$ in $H^1$. \item there
exists a sequence $y_n$ such that $u_n=u^1(x+y_n)+w_n$  with $w_n
\rightarrow 0$ in $H^1$.
\end{enumerate}
We have that $u^0, u^1 \in S$. In both cases Theorem \ref{stability}
holds. $\meproof$

\vspace{2mm} We give two  examples of functions $F$ which satisfy
the assumption (\ref{F0}), (\ref{F1}), (\ref{Fp'}):
\begin{enumerate}
\item $F(s)=-\frac{1}{4}|s|^4+\frac{1}{5}|s|^5$\ with $N=3$,\\
\item $F(s)=-\frac{|s|^q}{1+|s|^{q-p}}$\ with $2<p<
q<2+\frac{4}{N}$.
\end{enumerate}

With this potential we have that problem (\ref{eq}) has a solution
and that the  nonlinear Schr\"odinger equation

\begin{equation}\tag{\ref{schr}}
\left\{
\begin{array}{ll}
i\frac{\partial \psi}{\partial t}+\Delta\psi-F'(\psi)=0&(t,x)\in \R^+\times\R^N\\
&\\
\psi(0,x)=\psi_0(x)&
\end{array}
\right.
\end{equation}
admits a stationary solution which is orbitally stable.

\section{Appendix}

 \vspace{2mm}
 \noindent{\it Proof of the splitting
lemma.} We do it by steps.

\vspace{2mm} \noindent{\it Step I.} There exists $u^0\in H^1$ such
that $u_n\rightharpoonup u^0$ in $H^1$ and $u_0$ is a weak solution
of (\ref{eq}).

In fact, we have that $u_n$ is bounded in $L^2$ by hypothesis.
Furthermore, using (\ref{F0}), we also have that $u_n$ is bounded in
$H^1$. So there exist $u^0$ in $H^1$ such that $u_n\rightharpoonup
u^0$.

Now, because $u_n$ is a P-S sequence for $J_\lambda$, we have that
for all $\varphi\in C^\infty_0(\R^N)$,
\begin{equation}
\int \nabla u_n\varphi +\int F'(u_n)\varphi-\lambda \int
u_n\varphi\rightarrow 0.
\end{equation}
We have that, for any compact set $B$, $u_n\rightarrow u^0$ strongly
in $L^p(B)$ for $2\leq p<2^*$. Thus using (\ref{Fp'}) and the fact
that $u_n\rightharpoonup u^0$, we can conclude that $u^0$ is a weak
solution of (\ref{eq}).

\vspace{2mm} \noindent{\it Step II.} Setting $\psi_n=u_n-u^0$, we
have that
\begin{eqnarray}
||\psi_n||_{H^1}^2&=&||u_n||_{H^1}^2-||u^0||_{H^1}^2+o(1)\label{st2.1}\\
J(\psi_n)&=&J(u_n)-J(u^0)+o(1).\label{st2.2}\label{st2.2}
\end{eqnarray}
We have that $\psi_n\rightharpoonup0$ in $H^1$. Thus, obviously,
\begin{eqnarray*}
\int |u_n|^2&=&\int (u^0+\psi_n)^2=\int |u^0|^2+\int |\psi_n|^2+2\int u^0\psi_n=\\
&=&\int |u^0|^2+\int |\psi_n|^2+o(1)
\end{eqnarray*}
In the same way, we can proceed with $\displaystyle \int|\nabla
u_n|^2$, obtaining (\ref{st2.1})

To obtain (\ref{st2.2}), we prove that
\begin{equation}
\int F(u_n)-F(u^0)-F(\psi_n)\rightarrow0.
\end{equation}
For all $R>0$, we can write this integral as follows:
\begin{eqnarray*}
\int\limits_{\R^N}
F(u_n)-F(u^0)-F(\psi_n)&=&\int\limits_{B_R}[F(u^0+\psi_n)-F(u^0)]
-\int\limits_{B_R^C}F(u^0)+
\\
&&+\int\limits_{B_R^C}[F(u^0+\psi_n)-F(\psi_n)]-\int\limits_{B_R}F(\psi_n).
\end{eqnarray*}
On every compact set, we have that $\psi_n\rightarrow0$ in $L^p$
for all $2\leq p<2^*$ and that, by (\ref{Fp'}),
\begin{equation*}
\int\limits_{B_R}[F(u^0+\psi_n)-F(u^0)]-\int\limits_{B_R}F(\psi_n)\rightarrow0
\text{ when }n\rightarrow0.
\end{equation*}
Easily, we have also that
\begin{equation*}
\int\limits_{B_R^C}F(u_0)\rightarrow0 \text{ when
}R\rightarrow\infty.
\end{equation*}
Finally, for some $0<\theta<1$,
\begin{equation*}
\int\limits_{B_R^C}[F(u^0+\psi_n)-F(\psi_n)]=\int\limits_{B_R^C}F'(\theta
u^0+\psi_n)u^0,
\end{equation*}
with $||u^0||_{L^p(B_R^C)}\rightarrow0$ strongly in $L^p(B_R^C)$
when $R\rightarrow\infty$ and $\theta u^0+\psi_n$ is bounded in
$L^p(B_R^C)$. By (\ref{Fp'}), we have also that this term vanishes
when $R$ is sufficiently large, and this proves Step II.

\vspace{2mm} \noindent{\it Step III.} Set $\psi_n=u_n-u^0$. If
$\psi_n\nrightarrow0$ in $H^1$ then there exists a sequence of
points $y_n\in\R^N$, with $|y_n|\rightarrow\infty$, and a function
$u^1\in H^1$, $u^1\neq 0$, such that
\begin{equation}
\psi_n(x+y_n)\rightharpoonup u^1\  \in H^1.
\end{equation}

Notice that, if $\psi_n\rightarrow 0$ strongly in $H^1$ the
splitting lemma is proved. Otherwise, we start to prove that, when
$n\rightarrow\infty$,
\begin{equation}
\int\limits_{\R^N}F'(u_n)u_n-F'(u^0)u^0-F'(\psi_n)\psi_n\rightarrow0.
\end{equation}
As usual, for a fixed  $R>0$, we have that both $u_n\rightarrow u^0$
and $\psi_n\rightarrow0$ in $L^p(B_R)$ as $n\rightarrow\infty$. So,
by (\ref{Fp'})
\begin{equation}
\int\limits_{B_R}F'(u_n)u_n-F'(u^0)u^0-F'(\psi_n)\psi_n\rightarrow0.
\end{equation}
Moreover, there exist $\theta,\eta$, $0<\theta,\eta<1$ such that
\begin{eqnarray*}
&&\int\limits_{B_R^C}F'(u_n)u_n-F'(u^0)u^0-F'(\psi_n)\psi_n=\\
&&=\int\limits_{B_R^C}[F'(u^0+\psi_n)-F(u^0)]u^0+\int\limits_{B_R^C}[F'(u^0+\psi_n)-F'(\psi_n)]\psi_n=\\
&&=\int\limits_{B_R^C}[F''(u^0+\theta\psi_n)-F''(\eta
u^0+\psi_n)]u^0\psi_n
\end{eqnarray*}
and we can conclude as above that for $R$ sufficiently large this
term vanishes.

Using that $u_n$ is a P-S sequence and that $u^0$ is a weak solution
of (\ref{eq}), we have that
\begin{eqnarray*}
||\psi_n||^2_{H^1}&=&\int |\nabla u_n|^2-\int |\nabla u^0|^2+\int |u_n|^2-\int |u^0|^2+o(1)=\\
&=&-\int [F'(u_n)u_n-F'(u^0)u^0]+\\&+&(\lambda+1)\int[|u_n|^2- |u^0|^2]+o(1)=\\
&=&-\int F'(\psi_n)\psi_n+(\lambda+1)\int |\psi_n|^2+o(1).
\end{eqnarray*}
Now, for a fixed $L>0$, we decompose $\R^N$ into a numerable union
of $N$-dimensional hypercubes $Q_i$, having edge $L$.

By (\ref{Fp'}), we have
$$||\psi_n||_{H^1}^2+o(1)=-\int F'(\psi_n)\psi_n+(\lambda+1)\int
|\psi_n|^2\leq$$ $$\leq C_1\sum_i\left[
||\psi_n||_{L^q(Q_i)}^q+||\psi_n||_{L^p(Q_i)}^p+
(\lambda+1)||\psi_n||_{L^2(Q_i)}^2\right] $$ $$\leq  C_1\sum_i\Big[
L^{N\left(\frac{p-q}{p}\right)}||\psi_n||_{L^p(Q_i)}^q+
||\psi_n||_{L^p(Q_i)}^p
+(\lambda+1)L^{N\left(\frac{p-2}{p}\right)}||\psi_n||_{L^p(Q_i)}^2\Big]$$
$$\leq C_2||\psi_n||_{H^1}^2\left[L^N\left(L^{\left(\frac{p-q}{p}\right)}||\psi_n||_{L^p(\R^N)}^{q-2}+
(\lambda+1)L^{\left(\frac{p-2}{p}\right)}\right)+d_n\right]$$
where $d_n=\sup_i||\psi_n||_{L^p(Q_i)}^{p-2}$ and $C_1,C_2$ are
positive constants. We can choose $L$ small enough such that for a
suitable $C_3$ we have
\begin{equation*}
C_3||\psi_n||_{H^1}^2\leq d_n||\psi_n||^2_{H^1}+o(1).
\end{equation*}
Because $\psi_n\nrightarrow 0$ in $H^1$, we must have $\inf d_n >
0$. Thus there exists an $\alpha>0$ and a sequence of index $i_n$
such that
\begin{equation}\label{st3}
||\psi_n||_{L^p(Q_{i_n})}>\alpha \text{ for all }n.
\end{equation}
We call $y_n$ the center of the hypercube $Q_{i_n}$. Because
$\psi_n\rightarrow0$ in $L^p(B)$ for any compact set $B$, we have
that $|y_n|\rightarrow \infty$ when $n\rightarrow\infty$.

Finally, we know that there exists a $u^1$ in $H^1$ such that
\begin{equation}
\psi_n(\cdot+y_n)\rightharpoonup u^1
\end{equation}
weakly in $u^1$ because the sequence $\psi_n(\cdot+y_n)$ is bounded
in $H^1$, and by (\ref{st3}) we conclude that $u^1\neq0$.

\vspace{2mm} \noindent{\it Step IV.} The function $u^1$ is a weak
solution of (\ref{eq})

We know that $\psi_n(x+y_n)\rightarrow u^1$ weakly in $H^1$ and
strongly in $L^p(B)$ for all $B$ compact, $2\leq p<2^*$. So it is
sufficient to prove that, for all $\varphi\in C^\infty_0(\R^N)$,
\begin{equation}
\int\limits_{\R^N} \nabla
\psi_n(x+y_n)\nabla\varphi(x)+F'(\psi_n(x+y_n))\varphi(x)
-\lambda\psi_n(x+y_n)\varphi(x)dx \rightarrow0.
\end{equation}
After a change of variables we obtain
\begin{equation*}
\int\limits_{\R^N} \nabla
\psi_n(x)\nabla\varphi(x-y_n)+F'(\psi_n(x))\varphi(x-y_n)
-\lambda\psi_n(x)\varphi(x-y_n)dx,
\end{equation*}
and, using that $u^0$ is a weak solution of (\ref{eq}) and that
$u_n$ is a P-S sequence, we have that
\begin{eqnarray*}
&&\int\limits_{\R^N} \nabla \psi_n(x)\nabla\varphi(x-y_n)
-\lambda\psi_n(x)\varphi(x-y_n)dx=\\
&&=-\int\limits_{\R^N} F'(u_n(x))\varphi(x-y_n)
+\int\limits_{\R^N} F'(u^0(x))\varphi(x-y_n)+o(1).
\end{eqnarray*}
Thus we prove that
\begin{equation}
\int\limits_{\R^N}[F'(\psi_n(x))-F'(u_n(x))+F'(u^0(x))]\varphi(x-y_n)\rightarrow0
\text{ as }n\rightarrow\infty.
\end{equation}
As usual, for fixed $B$, we can split this integral as follows
\begin{eqnarray*}
&&\int\limits_{\R^N}[F'(\psi_n(x))-F'(u_n(x))+F'(u^0(x))]\varphi(x-y_n)=\\
&&\int\limits_{B_R}[F'(u^0(x))-F'(u_n(x))]\varphi(x-y_n)+\int\limits_{B_R^C}F'(u^0(x))\varphi(x-y_n)+\\
&&\int\limits_{B_R^C}[F'(\psi_n(x))-F'(u_n(x))]\varphi(x-y_n)+
\int\limits_{B_R}F'(\psi_n(x))\varphi(x-y_n).
\end{eqnarray*}
All the integrals over $B_R$ are definitively 0 because $\varphi$
has compact support and $|y_n|\rightarrow\infty$. Moreover, we
observe that $\displaystyle
\int\limits_{B_R^C}|F'(u^0)|^{\frac{p}{p-1}}\rightarrow 0$ as
$R\rightarrow\infty$.

Finally, for some $0<\theta<1$,
\begin{eqnarray*}
&&\int\limits_{B_R^C}[F'(\psi_n(x))-F'(u_n(x))]\varphi(x-y_n)=\\
&&=-\int\limits_{B_R^C}F''(\psi_n(x)+\theta
u^0(x))\varphi(x-y_n)u^0(x)\rightarrow0,
\end{eqnarray*}
as usual.

\vspace{2mm} \noindent{\it Step V.} Conclusion.

We can now iterate this procedure by defining a function
$$\psi^1_n(x)=\psi_n(x+y_n)-u^1(x).$$ We have that
$||u_n||_{H^1}^2=||u^0||_{H^1}^2+||u^1||_{H^1}^2+||\psi^1_n||$. If
$\psi^1_n\rightarrow 0$ strongly in $H^1$, the lemma is proved.
Otherwise, we have that $\psi^1_n\rightharpoonup 0$ in $H^1$ and
there exist a sequence of point $y^1_n$ with
$|y^1_n|\rightarrow\infty$ and a function $u^2$ in $H^1$ such that
$\psi^1_n(x+y^1_n)\rightharpoonup u^2 $ in $H^1$. Furthermore $u^2$
is a weak solution of (\ref{eq}), and so on.

We can have a finite number of iterative steps. Indeed, there
exists an $\alpha>0$ such that
\begin{equation}\label{palle}
||u^j||_{H^1}>\alpha \text{ for all }j.
\end{equation}

Hence, by (\ref{st2.1}) and (\ref{st2.2}) we get the claim. Now we
prove (\ref{palle}). We know that every $u^j$ is a weak solution
of (\ref{eq}), so it belongs to the set
\begin{equation}
{\cal N}:=\left\{u\in H^1,\ u\neq 0\ :\ \int |\nabla u|^2+\int
F'(u)u-\lambda \int u^2=0 \right\}.
\end{equation}
We want to prove that
\begin{equation}\label{st5}
\inf_{u\in{\cal N}}||u||_{H^1}=\alpha>0
\end{equation}

Notice that, because $\lambda<0$, then we can endow $H^1$ with the
following equivalent norm:
\begin{equation}
|||u|||_{H^1}=\int|\nabla u|^2-\lambda\int| u|^2.
\end{equation}
We suppose, by contradiction, that there exists a sequence
$w_n\in{\cal N}$ with $|||w_n|||_{H^1}\rightarrow0$. We can set
$w_n=t_nv_n$ with $|||v_n|||_{H^1}=1$, thus $t_n\rightarrow0$. We
have
\begin{eqnarray*}
0&=&\int |\nabla w_n|^2+\int F'(w_n)w_n-\lambda \int w_n^2=|||w_n|||_{H^1}+\int F'(w_n)w_n=\\
&=&t_n^2+t_n\int F'(w_n)v_n.
\end{eqnarray*}
Thus,
\begin{eqnarray*}
t_n&=&-\int F'(w_n)v_n\leq c_1 \int |t_nv_n|^{p-1}v_n+c_2 \int |t_nv_n|^{q-1}v_n\leq\\
&&\leq c_1 t_n^{p-1}\int |v_n|^p+c_2 t_n^{q-1} \int |v_n|^q;\\
1&\leq& c_1 t_n^{p-2}\int |v_n|^p+c_2 t_n^{q-2} \int |v_n|^q,
\end{eqnarray*}
and this lead to a contradiction. Indeed, if $v_n$ is bounded in
$H^1$ then it is bounded in $L^p$ for all $2\leq p<2^*$, and by
hypothesis $t_n\rightarrow0$. $\meproof$

\end{document}